\documentclass [12pt]{article}
\usepackage{graphicx,amssymb,amsfonts,latexsym,amsmath,amsthm,times}
\usepackage{epsfig}
\usepackage{fancyhdr}
\usepackage{color}
\usepackage[english]{babel}
\numberwithin{equation}{section}

\newtheorem{theorem}{Theorem}

\newtheorem{lemma}[theorem]{Lemma}

\DeclareMathOperator{\Arg}{Arg}
\DeclareMathOperator{\im}{Im}
\DeclareMathOperator{\re}{Re}
\DeclareMathOperator{\dr}{Dr}

\newcommand{\RR}{\mathbb{R}}
\newcommand{\CC}{\mathbb{C}}

\newcommand{\vphi}{\varphi}
\newcommand{\om}{\omega}
\newcommand{\Om}{\Omega}
\newcommand{\ep}{\varepsilon}

\newcommand{\pa}{\partial}
\newcommand{\ii}{\mathrm{i}}
\newcommand{\dd}{\mathrm{d}}

\newcommand{\ChiCpl}{\chi_{_{\im z >0} }}
\newcommand{\Adm}{Adm}

\newcommand{\wt}{\widetilde}

\begin{document}

%*******************************************************************
%TITLE OF THE JOURNAL - NAME OF THE ISSUE AS WELL AS YEAR AND PAGE NUMBERS
%*******************************************************************
%\footnotesize {\flushleft \mbox{\bf \textit{Math. Model. Nat. Phenom.}}}  \\
%\mbox{\textit{{\bf Vol. ??, No. ?, ????, pp. ?-?}}}

\thispagestyle{plain}

%*******************************************************************
%TITLE OF THE ARTICLE
%*******************************************************************
\vspace*{2cm} \normalsize \centerline{\Large \bf Nonlinear eigenvalue problem for optimal resonances in optical cavities}

\vspace*{1cm}
%*******************************************************************
%AUTHORS - THE CORRESPONDING AUTHOR NEEDS TO SPECIFY HIS/HER E-MAIL ADDRESS AS A FOOTNOTE
%*******************************************************************

\centerline{\bf I. M. Karabash$^a$ \footnote{Corresponding author.
E-mails: i.m.karabash@gmail.com, karabashi@mail.ru}}

\vspace*{0.5cm}

%*******************************************************************
%ADDRESS OF THE AUTHORS
%*******************************************************************
\centerline{$^a$ Institute of Applied Mathematics and Mechanics of
NAS of Ukraine,}
\centerline{R. Luxemburg str. 74, Donetsk 83114, Ukraine}

%*******************************************************************
%ABSTRACT
%*******************************************************************

\vspace*{1cm}

\noindent {\bf Abstract.} The paper is devoted to optimization of resonances in a 1-D open
optical cavity.
The cavity's structure is represented by its dielectric permittivity function $\ep (s)$.
It is assumed that $\ep (s)$ takes values in the range $1 \leq \ep_1 \leq \ep (s) \leq \ep_2$.
The problem is to design, for a given (real) frequency $\alpha$, a cavity having a resonance with
the minimal possible decay rate.
Restricting ourselves to resonances of a given frequency $\alpha$,
we define cavities and resonant modes with locally extremal
decay rate, and  then study their properties. We show that such locally extremal
cavities are 1-D photonic crystals consisting of alternating layers of two materials with extreme
allowed dielectric permittivities $\ep_1$ and $\ep_2$.
To find thicknesses of these layers, a nonlinear eigenvalue problem for locally extremal resonant modes is derived.
It occurs that coordinates of interface planes between the layers can be expressed via $\arg$-function
of corresponding modes. As a result, the question of minimization of the decay rate
is reduced to a four-dimensional problem of finding the zeroes of a function of two variables.

\vspace*{0.5cm}

%*******************************************************************
%KEYWORDS
%*******************************************************************

%*******************************************************************
%AMS SUBJECT CLASSIFICATION
%*******************************************************************

%*******************************************************************
%ARTICLE BODY
%*******************************************************************

\vspace*{1cm}

%*******************************************************************
%DO NOT FORGET TO RESET THE EQUATION COUNTER TO 0 AT THE HEAD OF EACH SECTION
%*******************************************************************
\setcounter{equation}{0}

\vspace*{0.5cm} \setcounter{equation}{0}
\section{Introduction}
\label{s intro}

A leaky optical cavity (or resonator) is a region of space within which the
electromagnetic field is well confined, but not completely
confined. Because of the leakage, each standing electromagnetic wave in the cavity decays exponentially
in time. The rate of energy loss and the frequency of these eigenoscillations can be characterized by
the corresponding complex eigenvalue $\om$, which is called a quasi-normal (QN) eigenvalue
or a resonance (the latter name is used in Quantum Mechanics and, sometimes, Engineering).
The imaginary part $\im \om$ of the QN eigenvalue corresponds to the decay rate of the standing wave,
the real part $\re \om$ to the frequency of oscillations.

Optical cavities with small and high decay rate are required for various
applications in optics \cite{V03} (including spontaneous emission control,
low-threshold lasers, low-power switches, and adiabatic wavelength conversion) because they generally enhance
intrinsically small light-matter interactions. Since light is essentially difficult to localize, it is
hard to realize small-sized optical cavities with strong light confinement. Recently,
however, the rapid progress has been made on this issue, see \cite{AASN03,NKT08} and references therein.
As a result, numerical aspects of emerging optimization problems for QN eigenvalues
have attracted great current interest among  specialists in Applied Mathematics \cite{HBKW08,KS08,BPSchZsch11}.

An attempt to provide an analytical background for QN eigenvalue optimization problems
lead the author to the development of the two-parameter perturbation
approach \cite{Ka12_pr}.
The paper \cite{Ka12_pr} is focused on mathematical details of the proofs, and
uses techniques of nonlinear functional analysis and multi-dimensional complex analysis.
The present communication is a continuation of the paper \cite{Ka12_pr}.
Its goal is to apply the two-parameter perturbation method to a concrete engineering problem
and, using this example, to make the main ideas understandable for specialists in applied sciences.
We consider the problem of optimization of a single QN eigenvalue in one-sided leaky optical cavity.
While the corresponding eigenvalue problem is derived under several simplifying assumptions,
it is proved to be very useful and have been extensively studied in Physics and numerical papers,
see e.g. \cite{U75,LLTY94,HBKW08}. We give rigorous mathematical definitions of locally extremal QN eigenvalues
and of optimal cavities. In Section \ref{ss def extr}, the existence of such cavities is shown.
Then the two-parameter perturbation method is introduced in Section \ref{s pert} and applied to derive a nonlinear equation for QN modes with
extremal properties in Section \ref{s nonlin}. We use this nonlinear equation to reduce the cavity's optimization,
which is essentially an infinite-dimensional problem,
to the four-dimensional problem of finding the zeroes of a specially introduced function $W$,
which depends only on one real and one complex variable.
The main conclusions are summarized in the last section.

Some of purely mathematical details of the presented method are sketched to make the paper more accessible for non-mathematicians.
However, a reader interested in rigorous analytical proofs will be able to recover them without difficulty
using the proofs in \cite{Ka12_pr}.

Note that, in Engineering and Physics literature,
the rate of energy loss is often characterized by other parameters, like Q-factor or lifetime of a QN mode,
which in many cases can be approximately expressed in terms of $\im \om$ and $\re \om$ \cite{LL84}.
Keeping this in mind, it is easy to rewrite the results of the paper in terms of these,
to some extend equivalent, parameters.

\vspace*{0.5cm} \setcounter{equation}{0}
\section{Physical settings and the problem statement} \label{s ph set}

Consider an isotropic non-dispersive transparent medium which
has the relative magnetic permeability equal to 1 everywhere,
but which is electrically inhomogeneous in one direction. The
electromagnetic structure of this medium is described by the
\emph{(dielectric) permittivity} $\ep$, which is $\ge 1$ and varies
only in the $x_3$-direction $\ep=\ep(x_3)$. When $\ep$ is piecewise
constant, this model is called a \emph{1-D photonic crystal}.

Under additional simplifying assumptions that electromagnetic waves
pass normally and the electric field $\mathbf{E}$ is in the
$x_2$-direction $\mathbf{E} = \{ 0, E, 0 \} $, the study of
electromagnetic wave propagation can be reduced to the scalar 1-D
wave equation
\begin{equation} \label{e weE}
\pa_{x_3}^2 E (x_3,t) = \frac{\ep (x_3)}{c^2} \pa_t^2 E (x_3,t) ,
\end{equation}
where the $x_2$-component $E$ of the electric field depends only on
$x_3$, and $c$ is the speed of light in vacuum.

To introduce a one-sided open optical cavity, assume that the medium
has the above structure in the slab $0<x_3<l$  with piecewise continuous
$\ep$ (that is, for a certain partition $0=s_0<s_1 < \dots < s_{n-1} < s_n=l$,
$\ep (s)$ is continuous in the intervals $s_{j-1} < s < s_{j}$
and has one-sided limits at points $s_{j}$).
The outer medium is assumed to be  a perfect conductor for $x_3<0$, and
a homogeneous dielectric with constant permittivity $\ep_\infty \ge 1$ for $x_3>l$.
The boundary condition for the electric field on the
plane $x_3=0$ is $E = 0$.  Since the medium is homogeneous for
$x_3>l$, the waves radiated from the cavity are never reflected
back. If all sources of waves are inside the cavity $0<x_3<l$, the
waves satisfy the radiation boundary condition
$\pa_{x_3} E + \frac{\sqrt{\ep_\infty}}{c} \, \pa_t E = 0$ at the interface $x_3 = l$.
For \emph{a monochromatic field}
$E (x_3,t) = E (x_3) e^{-\ii \om t}$ in the cavity, these settings lead to the
eigenvalue problem
\begin{equation} \label{e epE}
E'' (x_3) = - \om^2 \frac{\ep(x_3)}{c^2} E (x_3)
\end{equation}
equipped with the boundary conditions
\begin{equation}
E (0)  =  0, \ \ \ \ii \om \frac{\sqrt{\ep_\infty}}{c} E (l) =
E'(l) . \label{e bcE}
\end{equation}

Complex eigen-parameters $\om=\alpha - \ii \beta$ are called
\emph{quasi-normal (QN) eigenvalues}. The real part $\alpha=\re \om
$ is the \emph{frequency} of the standing wave $E (x_3)
e^{-\ii \om t}$. Since the energy leaks from the cavity, the
oscillations are decaying. The minus imaginary part $\beta = -
\im \om$ is positive for all QN eigenvalues and is called
the \emph{decay rate}.

%\subsection{The problem statement} \label{ss problem}

Let $1 \leq \ep_1 < \ep_2$.
Restricting ourselves to cavities with $\ep$ in the range
\begin{equation} \label{e <ep<}
 \ep_1 \leq \ep (x_3) \leq \ep_2,
\end{equation}
consider the following optimization problem: \emph{devise, for a given frequency
range $\alpha_1 \leq \alpha \leq \alpha_2$, a cavity that has a QN
eigenvalue $\om = \alpha - \ii \beta$ with the minimal possible decay
rate $\beta$}.

\vspace*{0.5cm} \setcounter{equation}{0}
\section{QN eigenvalues: basic properties and extrema's definitions}
\label{s math set}

\subsection{QN eigenvalues and their multiplicities}

Since the QN eigenvalue problem is one-dimensional, from now on we
will write $s$ instead of $x_3$. Let the length $l$ of the cavity
and the permittivity $\ep_\infty $ of the outer medium $s>l$ be
fixed. Then the cavity is completely described by the permittivity
$\ep(s)$ in the interval $0 \le s \le l$ and, for brevity's sake, we will speak on
\emph{QN eigenvalues of the cavity} $\ep$.

By $\vphi(s) = \vphi (s,z;\ep)$ and $\psi(s) = \psi (s,z;\ep)$ we denote the solutions of
$E'' (s) = - z^2 \frac{\ep(s)}{c^2} E (s)$ satisfying
\begin{equation*} \label{e phi psi}
\vphi (0, z; \ep ) = \pa_s \psi (0, z; \ep) = 1 , \ \ \pa_s \vphi
(0, z; \ep ) = \psi (0, z; \ep) = 0 .
\end{equation*}
In particular, $\psi $ satisfies the integral equation
\[
\psi (s)  =  s - \frac{\om^2}{c^2} \ \int_0^s (s-\tau) \ \ep (\tau) \
\psi (\tau) \ \dd \tau , \quad 0 \leq s \leq l .
\]

The QN eigenvalues of $\ep$ are exactly the roots of the equation
\begin{equation} \label{e F}
F(z) = 0, \ \ \ \text{where } F(z)  = F (z;\ep) =
\ii z \frac{\sqrt{\ep_\infty}}{c} \psi (l,z;\ep) -  \pa_s \psi (l,z;\ep) .
\end{equation}
The function $F$ is analytic in $z$ in the whole complex plane $\CC$.
A nontrivial solution $E (s)$ to the problem (\ref{e epE}), (\ref{e bcE}) is called
\emph{a QN mode} corresponding to a QN eigenvalue $\om$. Recall that a solution $E$ is called trivial if
$E(s) = 0$ for all (more precisely, almost all) $s$ in the interval $0 \le s \le l$.

All QN modes of $\ep$ have the form $C \psi (s,\om; \ep)$, where $C$ is an arbitrary non-zero constant.
Hence the geometric multiplicity of any QN eigenvalue $\om$ equals 1. In
the following, \emph{the multiplicity of} $\om$ means its algebraic
multiplicity. By definition, the algebraic multiplicity of a QN
eigenvalue $\om$ is its multiplicity as a zero of $F$ \cite{Kel51}.
The $\om$-dependent boundary condition
$\ii \om \frac{\sqrt{\ep_\infty}}{c} E (l) =  E'(l)$ makes the QN eigenvalue problem
non-self-adjoint. In particular, QN eigenvalues may be
\emph{degenerate}, i.e., of multiplicity $\geq 2$, see e.g. \cite{KN79,GP97,Sh97}.
Concerning the contemporary theory of spectral problems with an eigen-parameter in boundary conditions, we refer to
\cite{AHM07} and references therein.

Concerning other properties of QN eigenvalues, it is easy to derive that:
(i) the multiplicity of each QN eigenvalue is finite,
(ii) the set of QN eigenvalues is symmetric with respect to (w.r.t.) imaginary axis
(together with multiplicities),
(iii) QN eigenvalues are isolated, $\infty$ is their only
possible accumulation point.
The property (ii) allows one to restrict the study of QN eigenvalues to the case $\re \om \ge 0$.

For a homogeneous cavity with constant permittivity $\ep \geq 1$, the QN eigenvalues are roots of
$ \tan (- \om \sqrt{\ep} \ l /c
) = \ii \sqrt{\ep/\ep_\infty} $.
When $\ep \neq \ep_\infty$, they are non-degenerate and form a uniformly spaced
sequence
\begin{equation} \label{e om n}
\om_n = - \ii \frac{c}{2l\sqrt{\ep}}
 \ln \left| \frac{\ep + \ep_\infty}{\ep - \ep_\infty } \right| + \frac{\pi c}{l\sqrt{\ep}}
 \left\{
\begin{array}{ll}
n , & \text{if } \ep < \ep_\infty\\
n+1/2, & \text{if } \ep > \ep_\infty
\end{array} \right. , \ \ \ n = \dots, -1,0,1, \dots.
\end{equation}
The case when the space outside the cavity is vacuum (i.e., $1 = \ep_\infty \leq \ep$) is discussed in detail in
\cite{U75}.
In the case $\ep = \ep_\infty$,
\emph{there are no QN eigenvalues} since the medium is homogeneous for $s>0$
(the energy of the initial disturbance
localized in the region $0<s<l$ escapes this region in finite time, see e.g. \cite{CZ95} for detailed explanations).

\subsection{Local and global extrema for admissible frequencies}
\label{ss def extr}

\emph{The family $\Adm$ of admissible cavities} is defined by
condition (\ref{e <ep<}). Mathematically, it is convenient to
include in $\Adm$ all $L^\infty$ functions satisfying (\ref{e
<ep<}). \emph{An admissible direction} $h$ for an admissible cavity
$\ep$ is a function on the interval $0 \le s \le 1$ such that the cavity $\ep + h $ is admissible
(we call such $h$ a direction since, for all constants $C$ in the range $0 < C \le 1$,
the cavity $\ep + C h $ is admissible).
We say that $\om$ is \emph{an admissible QN eigenvalue} if it is a
QN eigenvalue of some admissible cavity, and say that $\alpha$ is
\emph{an admissible frequency} if $\alpha = \re \om$ for some
admissible QN eigenvalue $\om$.

To get some understanding of admissible frequencies, consider homogeneous cavities.
Changing constant $\ep$ in the range $\ep_1 \leq \ep \leq \ep_2$ and taking $\om_n$ given
by (\ref{e om n}), we see that the frequencies $\re \om_n$ are admissible.
They form a sequence of intervals. For $n$ large enough, the interval produced by $\re \om_n$
overlaps with the next one produced by $\re \om_{n+1}$.
Hence frequencies $\alpha$ with $|\alpha|$  large enough are admissible.
For instance, let us look closer at the case $\ep_\infty \leq \ep_1 $.
When constant $\ep$ runs over the interval $\ep_1 < \ep \leq \ep_2$, the real part of $\frac{l \om_n }{\pi c}$ runs over the interval
$(n+1/2) \ep_2^{-1/2} \leq \alpha < (n+1/2) \ep_1^{-1/2}$. When $n \geq \frac{1}{(\ep_2/\ep_1)^{1/2} -1} - \frac 12$, there
is no gap between the $n$-th and $(n+1)$-st intervals. Hence a frequency $\alpha$ is admissible whenever
\[
|\alpha| \geq \frac{\pi c}{l \sqrt{\ep_2}} \left(\frac 12 +  \left\lceil \frac{1}{(\ep_2/\ep_1)^{1/2} -1} - \frac 12 \right\rceil  \right) ,
\]
(here $\lceil s \rceil $ is the smallest integer that is greater or equal to $s$).

If an admissible cavity $\ep$ has a QN eigenvalue
$\om$, we say that $\{\om,\ep\}$ is \emph{an admissible pair}.
\emph{The decay rate} of an admissible pair is defined by
\[
\dr (\om,\ep) = - \im \om .
\]
While $\ep$ does not participate in the right hand side, it participates in the domain of definition
of the functional $\dr$. This is essential for the next definition.

We say that an admissible pair $\{\om,\ep\}$ is \emph{a local minimizer of $\dr$ for a frequency $\alpha$}
if, for all sufficiently small admissible directions $\delta \ep$ that
keep the frequency of a perturbed QN eigenvalue $\om +\delta \om$ of the cavity $\ep + \delta \ep$
equal to $\alpha$,
the decay rate of $\{\om + \delta \om, \ep + \delta \ep \}$ is $\geq \dr (\om,\ep)$.
The rigorous form of this condition is: there exist $\delta >0$ such that
$\dr (\wt \om , \wt \ep ) \geq \dr (\om,\ep)$ for any admissible
pair $\{ \wt \om , \wt \ep \}$ satisfying
\begin{eqnarray*}
| \wt \ep (s) - \ep (s)| & < \delta &\ \ \mbox{ for almost all } \ 0<s<l, \\
| \wt \om - \om| & < \delta , & \mbox{ and } \re \wt \om = \alpha.
\end{eqnarray*}
\emph{Local maximizers} of $\dr$ for a particular frequency are defined in a similar way.

We say that $\beta_{\min} (\alpha)$ is \emph{a minimal decay rate} and $\ep$ is
\emph{an optimal cavity} for a frequency $\alpha$, if the pair $\{\alpha - \ii \beta_{\min} (\alpha) , \ep \}$
is admissible and $\beta_{\min} (\alpha) \leq \dr (\wt \om, \wt \ep)$ for any admissible pair $\{ \wt \om, \wt \ep \}$
with $\re \wt \om = \alpha$.
In other words,
such a pair $\{\alpha - \ii \beta_{\min} (\alpha), \ep \}$ is a global minimizer (over the admissible family $\Adm$) of
the decay rate for a frequency $\alpha$ .
We will also say that the corresponding
\emph{QN eigenvalue $\alpha - \ii \beta_{\min} (\alpha) $ is optimal} for the frequency $\alpha$.

\vspace*{0.25cm}
\begin{lemma} \label{l exist}
For each admissible frequency, there exists an optimal cavity.
\end{lemma}

The existence of optimal cavities follows from the fact that
if admissible QN eigenvalues $\om_n $ approach a complex number $\om$, then
$\om$ is an admissible QN eigenvalue. To explain this fact, consider
admissible cavities $\ep^{(n)}$ that produce the QN eigenvalues $\om_n$.
Using the Gronwall-Bellman inequality and assumption (\ref{e <ep<}), it is easy to show that the sequence of related $\psi$-solutions
$\psi_n (s) = \psi (s,\om_n; \ep^{(n)})$ is uniformly bounded in $s$ and $n$ together with the sequences of
their $s$-derivatives $\psi_n '$ and $\psi_n ''$. Hence some subsequence $\psi_{n_j}$ converges uniformly to
a continuous function $\wt \psi (s)$ (due to the compact embedding of $W^{2,\infty}[0,1]$ into $C [0,1]$).
Then applying the Banach-Alaoglu theorem for $L^\infty (0,1)$ to the sequence $\ep^{(n_j)} $, one
can pass to a limit in the integral reformulation of (\ref{e epE})-(\ref{e bcE})
(see \cite{Ka12_pr}, and also \cite{Kr51} for compactness arguments applied to a self-adjoint eigenvalue optimization problem).

%such that subsequence  and $psi are
%we use the integral form of b.v.p. (\ref{e epE})-(\ref{e bcE})
%\begin{eqnarray}
%E (s) & = & E'(0) s - \frac{\om^2}{c^2} \ \int_0^s (s-s) \ \ep (s) \
%E (s) \ \d s , \quad 0 \leq s \leq 1,
%\label{e int ep} \\
%E (l) & - & \i \frac{\om}{c} \int_0^1 \ep (s) E (s) \d s = 0
% \ . \label{e int bc1}
%\end{eqnarray}

\vspace*{0.5cm} \setcounter{equation}{0}
\section{Two-parameter perturbations}
\label{s pert}

We study properties of \emph{local extremizers} (i.e., local minimizers and maximizers) via the
perturbation theory for
QN eigenvalues using the following arguments. Assume that
$\{\om,\ep\}$ is a local minimizer. If we consider all small
admissible directions $\delta \ep$ that keep $\re (\om + \delta \om)$ equal to
$\alpha$, then corresponding corrections $\delta \beta = - \im \delta
\om$ of the decay rate must be positive. The key point is to show
that there are enough various $\delta \ep$ with the above property to
impose strong restrictions on extremal $\ep$.

For non-degenerate QN eigenvalues, the first order
perturbation correction can be calculated by the standard sensitivity analysis.
However our method requires
a more delicate approach due to following
reasons: (i) QN eigenvalues may be degenerate, (ii) we are
interested only in the perturbed QN eigenvalues that stay on the
line $\re \om = \alpha$.
Indeed, even if one is able to find the first order correction that
moves a QN eigenvalue along the line $\re \om = \alpha$ using singular
perturbation theory (see e.g. \cite{MBO97} and references therein), there are
no any guarantee that higher order corrections does not contribute
to the real part of $\om$.

To resolve these difficulties we study two-parameter perturbations of QN eigenvalues
considering them as $z$-roots of the equation $F (z,\ep+ \delta \ep)=0$.
This requires computation of the (complex) partial derivative $\pa_z F (z;\ep)$ and
the directional derivatives $\pa_\zeta F (z; \ep + \zeta h)$ for admissible pairs $\{\om, \ep \}$.
Here and below $\zeta$ is a complex number, and $h$ is an arbitrary (measurable and essentially bounded) function on
$0< s < l$, which is called \emph{a direction}.

The derivative $\pa_\zeta F (z; \ep + \zeta h)$ exists for an arbitrary direction.
More precisely, the Maclaurin series for the solution $\psi$
(see e.g. \cite{KK68_II}) can be used to get the Maclaurin series for
$F(z ; \ep)$ and, in turn, to show that for arbitrary $n$ directions $h_1 (s)$, \dots, $h_n (s)$,
the function
%$Q (z,\zeta_1,\zeta_2,\dots,\zeta_n) = F (\om + z; \ep + \zeta_1 h_1  + \zeta_2 h_2 + \dots + \zeta_n h_n)$ is an
$ F (\om + z; \ep + \zeta_1 h_1  + \zeta_2 h_2 + \dots + \zeta_n h_n)$ is an
analytic function of n+1 complex variables $z$, $\zeta_1$, \dots, $\zeta_n$
(moreover, the map $\{z,\ep\} \mapsto F (z ;\ep)$ is analytic on $\CC \times L^\infty (0,1)$).

Differentiating the integral equation for $\psi$ and solving the produced differential equations,
it is possible to express  $\pa_z F (z;\ep)$ and $\pa_\zeta F (z; \ep + \zeta h)$ via
the fundamental solutions $\vphi$ and $\psi$.
%Obviously,
%\begin{equation} \label{pa_z}
%\pa_z F (z;\ep) = \ii \om \sqrt{\ep_\infty} c^{-1} \pa_z \psi (1,z;\ep) - \pa_z \pa_s \psi (1,z;\ep).
%\end{equation}
%Differentiating the integral equation for $\psi$
%\[
%\psi (s,z;\ep) = s - \frac{z^2}{c^2} \int_0^s (s-s) \ep (s) \psi (s) \d s,
%\]
%we see that $\pa_z \psi (s,z;\ep)$ and $\pa_\zeta \psi (s,z;\ep+\zeta h)$ are solutions to
%the problem $y(0)=y'(0)=0$ for the equation
%$y''+ z^2 c^{-2} \ep y = f$ with the nonhomogeneous term
%$f(s)=-2 \om c^{-2} \ep (s)\psi(s)$ and $f(s) = - \om^2 c^{-2} h(s) \psi(s)$, respectively.
%Using these solutions to differentiate (\ref{e F}), one gets
\emph{When $z=\om$ is a QN eigenvalue of $\ep$, the obtained expressions can be simplified to}
\begin{eqnarray} \label{e pa zeta F}
\pa_\zeta F (\om; \ep + \zeta h)  & = & \frac{\om^2}{c^2 \psi (l,\om;\ep)} \int_0^l  \psi^2 (s,\om;\ep) h(s) \dd s , \\
\pa_z F (\om;\ep) & = & \frac{2 \om}{c^2 \psi (l,\om;\ep)} \int_0^l  \psi^2 (s,\om;\ep) \ep (s) \dd s +
\ii \frac{\sqrt{\ep_\infty}}{c} \psi (l,\om;\ep). \label{e pa z F}
\end{eqnarray}
The last equality yields the following lemma.
\vspace*{0.25cm}
\begin{lemma} A QN eigenvalue $\om$ of $\ep$ is degenerate exactly when
\[
\int_0^l  \psi^2 (s,\om;\ep) \ep (s) \dd s + \ii \frac{c\sqrt{\ep_\infty}}{2 \om} \psi^2 (l,\om;\ep) = 0 .
\]
\end{lemma}

Note that
\[
\psi (l,\om;\ep) = \left[ \ii \om \frac{\sqrt{\ep_\infty}}{c} \vphi (l,\om;\ep) - \vphi' (l,\om;\ep) \right]^{-1} \neq 0.
\]

Let us consider the first-order correction for $\om$ in the case when \emph{$\om$ is a non-degenerate QN eigenvalue  of $\ep$.}
Small perturbations of $\ep$ lead to small perturbations of $\om$ and the perturbed QN eigenvalues
remains non-degenerate.  Taking $\delta \ep = \zeta h$ with a fixed perturbation direction $h=h(s)$, it is easy to show
that the perturbed QN eigenvalue $\om + \delta \om$ is an analytic function $\Om (\zeta)$ of $\zeta$ in a vicinity of $\zeta_0=0$.
%(moreover, $\om + \delta \om$ is an analytic in $\delta \ep$ map defined in a vicinity of the $L^\infty (0,1)$-origin
%and taking values in $\CC$).
\emph{The approximation of $\Omega (\zeta)$ to the first-order term is given by}
\begin{equation} \label{e Om ndeg}
\Om (\zeta) = \om - \zeta  \frac{\om \int_0^l  \psi^2 (s) h(s) \dd s }
{2  \int_0^l  \psi^2 (s) \ep (s) \dd s + \ii \frac{c\sqrt{\ep_\infty}}{\om} \psi^2 (l)} + O (\zeta^2) .
\end{equation}
This formula, as well as higher order corrections in the non-degenerate case, are known \cite{LLTY94}.

Consider the case  when \emph{$\om$ is degenerate}. Perturbations of a multiple eigenvalue may lead to its
splitting in several eigenvalues. The general splitting picture is quite complicated and is described
by one or several series (Puiseux series) in fractional powers of $\zeta$ \cite{RS78IV,K80,MBO97}.
For our needs, it is enough to consider only the generic case when the direction $h$ satisfies
$\pa_\zeta F (\om; \ep + \zeta h) \neq 0$ . This condition is equivalent to
\begin{equation} \label{e pa h neq 0}
 \int_0^l  \psi^2 (s,\om;\ep) h(s) \dd s \neq 0 .
\end{equation}

Assume that $\om$ is a QN eigenvalue of $\ep$
of multiplicity $m$, and that perturbation's direction $h$ satisfies (\ref{e pa h neq 0}).
To study the equation $F (\om + z; \ep + \zeta h) = 0$, we introduce the function
$Q (z,\zeta) = F (\om + z; \ep + \zeta h)$ and write it as a
power series in $z$ and $\zeta$. The equation $Q=0$ takes the form
\[
z^m \frac{\pa_z^m Q (0,0)}{m!}  + \zeta \pa_\zeta Q (0,0)  + o (z^m) + o (\zeta) = 0
\]
for $z$ and $\zeta$ going to $0$. Note that $\pa_z^m Q (0,0) \neq 0$ since $\om$ has multiplicity $m$.
Considering $z$ as a function of $\zeta$, we see that
$|z(\zeta)| \leq C |\zeta|^{1/m}  $ for $\zeta$ small enough with some constant $C$. So
$o (z^m)$ is also $o (\zeta)$. Hence,
\[
z^m(\zeta) = - \zeta \frac{m! \pa_\zeta Q (0,0) }{\pa_z^m Q (0,0)} + o(\zeta) .
\]
According to (\ref{e pa zeta F}) and the assumption on $h$,
\[
\pa_\zeta Q (0,0) = \frac{\om^2}{c^2 \psi (l,\om;\ep)} \int_0^l  \psi^2 (s,\om;\ep) h(s) \dd s \neq 0 .
\]
Taking the $m$-th root,
we get a formal first-order correction formula for a perturbed QN eigenvalue $\Om (\zeta)$
%\begin{equation} \label{e Om mth}
%\Om (\zeta) = \om + \left( - \zeta \frac{m! \om^2  \int_0^l  \psi^2 (s,\om;\ep) h(s) \dd s }
%{c^2 \psi (l,\om;\ep) \pa_z^m F (\om,\ep)} \right)^{1/m} + o (|\zeta|^{1/m}) ,
%\end{equation}
\begin{equation} \label{e Om mth}
\Om (\zeta) = \om + \left( \zeta \ C_1 (\om,\ep,m)  \int_0^l  \psi^2 (s,\om;\ep) h(s) \dd s
 \right)^{1/m} + o (|\zeta|^{1/m}) ,
\end{equation}
where
\begin{equation} \label{e C1}
C_1 (\om,\ep,m) =  -  \frac{m! \om^2  }
{c^2 \psi (l,\om;\ep) \pa_z^m F (\om,\ep)}  .
\end{equation}
The $m$-th derivative $\pa_z^m F (\om,\ep)$ can also be expressed through the fundamental solutions $\vphi$ and $\psi$,
but
we are not going to use such a refinement in this paper.
Note that (\ref{e Om mth}) turns into (\ref{e Om ndeg}) for $m=1$.

\emph{Now consider two-parameter perturbations of $\om$,} assuming that the directions
$h_1$ and $h_2$ are such that $\pa_{\zeta_1} F (z; \ep + \zeta_1 h_1)$ and $\pa_{\zeta_2} F (z; \ep + \zeta_2 h_2)$
are nonzero and
\begin{equation} \label{e arg pa}
0 < \Arg
\pa_{\zeta_2} F (z; \ep + \zeta_2 h_2) - \Arg \pa_{\zeta_1} F (z; \ep + \zeta_1
h_1) < \pi
 \end{equation}
for a suitable branch $\Arg$ of the $\arg$-function (here and below $\arg z$ is the complex argument or phase of a nonzero complex number $z$).
Condition (\ref{e arg pa}) means that $ \pa_{\zeta_1} F (z; \ep + \zeta_1 h_1)$ and $ \pa_{\zeta_2} F (z; \ep + \zeta_2 h_2)$,  perceived
as vectors must not be oppositely directed and, additionally, must be ordered such that (\ref{e arg pa}) holds.
Considering a perturbed eigenvalue $ \om + z $ of $\ep+\zeta_1 h_1 +\zeta_2 h_2$ as a function of $\zeta_1$ and $\zeta_2$,
we write the equation
\[
Q (z , \zeta_1 , \zeta_2 ) = 0
\]
for the correction $z=z (\zeta_1,\zeta_2)$
defining $Q$ as
\[
Q (z,\zeta_1,\zeta_2) = F (\om + z, \ep + \zeta_1 h_1  + \zeta_2 h_2 ).
\]
Denote $\zeta = \{ \zeta_1 , \zeta_2 \}$, $|\zeta| = (\zeta_1^2 + \zeta_2^2)^{1/2}$, and $\mathbf{0} = \{ 0,0,0\} $.
When $|\zeta|$ is small, we have
\begin{equation} \label{e z m two}
z^m = - \zeta_1 \frac{m! \pa_{\zeta_1} Q (\mathbf{0})}{\pa_z^m Q \mathbf{0}}
- \zeta_2 \frac{m! \pa_{\zeta_2} Q (\mathbf{0}) }{\pa_z^m Q (\mathbf{0})}
+ o(z^m) + o (\zeta_1) + o (\zeta_2) .
\end{equation}

\emph{Assume additionally that}
\begin{equation} \label{e as z1 z2}
 \zeta_1 \text{ and } \zeta_2 \text{ are nonnegative and
of the same order}.
\end{equation}
That is, $\zeta_1 $ and $\zeta_2$ go to zero such
that their ratios $\zeta_1/\zeta_2$ and $\zeta_2/\zeta_1$ are
bounded. Under this condition, (\ref{e z m two}) and (\ref{e arg pa}) yield that
$|z (\zeta_1, \zeta_2)| \leq c |\zeta|^{1/m}$ for small $|\zeta|$.
Hence,
\begin{equation} \label{e Om m two}
z (\zeta_1,\zeta_2) = \left( - \zeta_1 \frac{m! \pa_{\zeta_1} Q (\mathbf{0})}{\pa_z^m Q (\mathbf{0})}
- \zeta_2 \frac{m! \pa_{\zeta_2} Q (\mathbf{0}) }{\pa_z^m Q (\mathbf{0})} \right)^{1/m} + o (|\zeta |^{1/m})  .
\end{equation}
Appending the values of partial derivatives
provided by (\ref{e pa zeta F}), one gets the first-order approximation for perturbed QN eigenvalues
%\begin{equation} \label{e Om mth two}
%\Om (\zeta_1,\zeta_2) = \om + \left( -  \frac{m! \om^2  \int_0^l  \psi^2 (s,\om;\ep) [\zeta_1 h_1 (s) + \zeta_2 h_2 (s)] \dd s }
%{c^2 \psi (l,\om;\ep) \pa_z^m F (\om,\ep) } \right)^{1/m} + o (|\zeta|^{1/m}) .
%\end{equation}
\begin{equation} \label{e Om mth two}
\Om (\zeta_1,\zeta_2) = \om + \left( C_1 (\om,\ep,m) \int_0^l  \psi^2 (s,\om;\ep) [\zeta_1 h_1 (s) + \zeta_2 h_2 (s)] \dd s
\right)^{1/m} + o (|\zeta|^{1/m})
\end{equation}
with $C_1$ defined by (\ref{e C1}).

For sufficiently small positive $\zeta_1$, $\zeta_2$ satisfying (\ref{e as z1 z2}),
the $m$ branches of $m$-th root produce $m$ different perturbed QN eigenvalues, each of which is non-degenerate.
The proof of this fact can be given along the lines suggested in the paper \cite{Ka12_pr}, where a slightly
different eigenvalue problem was considered.
Condition (\ref{e arg pa}) ensures that the coefficient under the $m$-th root is nonzero.
It is important for the sequel that this coefficient is a linear functional of perturbation directions $h_1$ and $h_2$.

Formula (\ref{e Om mth two}) will be applied to the study of the QN eigenvalue optimization
via the  lemma in the following section.

\vspace*{0.5cm} \setcounter{equation}{0}
\section{Admissible directional derivatives  and arguments of eigenvalue corrections}

Let us fix an admissible pair $\{\om,\ep\}$.
%Then a direction $h$ is called admissible for $\ep$ if the cavity
%$\ep + h$ is admissible.
The set of \emph{admissible directional derivatives} of $F$ at $\{ \om,\ep \}$ is defined as the set
formed by complex numbers $\pa_\zeta F (\om;\ep+\zeta h)$ with
$h$ running through the set of admissible directions (i.e., we take all $h$ such that $\ep +h$ is an admissible
cavity).
The set of admissible directional derivatives contains $0$ and is \emph{convex} (that is, if
$\xi_1$ and $\xi_2$ are admissible  directional derivatives, then $\gamma \xi_1 + (1-\gamma) \xi_2$ is
so for all $0 \le \gamma \le 1$). We say that the set of admissible directional derivatives \emph{contains
a neighborhood of zero} if
it contains all complex numbers $\xi$ with modulus $|\xi|$ less than a certain
positive number. Let $\Arg_0$ be the branch of the $\arg$-function that takes values
in the interval $-\pi <  \Arg_0 z \le \pi$.

\vspace*{0.25cm}
\begin{lemma} \label{l dir der}
Assume that the set of admissible directional derivatives of $F$ at an admissible pair $\{ \om,\ep \}$
contains a neighborhood of zero. Then for any $\theta$ in the range $-\pi <  \theta \le \pi$, there exists
a sequence of admissible pairs $\{ \om_n , \ep^{(n)}\}$ such that
\item[(i)] $\Arg_0 (\om_n - \om) = \theta$
for all $n=1,2,\dots$,
\item[(ii)] $\om_n $ tends to $\om$ and $\ep^{(n)}$ tends to $\ep$
(the latter means that the essential suprema of $|\ep (s) - \ep^{(n)} (s)|$
tends to zero as $n$ goes to $\infty$).
\end{lemma}

In particular, the lemma yields that \emph{if $\{ \om,\ep \}$ is a local extremizer of $\dr$ for
the frequency $\re \om$,
than the set of admissible directional derivatives of $F$ at $\{ \om,\ep \}$ cannot
contain a neighborhood of zero.}

To derive the lemma, note that its assumptions allows us to chose admissible directions
$h_1$ and $h_2$ and the branch of the m-th root in (\ref{e Om mth two})
such that
\begin{itemize}
\item[(A1)] assumption (\ref{e arg pa}) is valid for $h_1$ and $h_2$,
\item[(A2)]   for positive $\zeta$,
the corresponding first order corrections to $\om$,
\[
\left( \zeta C_1 (\om,\ep, m) \int_0^l  \psi^2 (s,\om;\ep) h_1 (s) \dd s \right)^{1/m}
 \text{and }
\left( \zeta C_1 (\om,\ep, m)    \int_0^l  \psi^2 (s,\om;\ep) h_2 (s) \dd s
\right)^{1/m} \!\!,
\]
have the arguments slightly smaller and, respectively, slightly greater than $\theta$,
\item[(A3)] the chosen branch of the m-th root is continuous in the sector of complex plane
covered by numbers
\[
C_1 (\om,\ep, m)  \int_0^l  \psi^2 (s,\om;\ep) [\zeta_1 h_1 (s) + \zeta_2 h_2 (s)] \dd s
\]
while $\zeta_1$ and $\zeta_2$ run through positive numbers.
\end{itemize}
This allows us to apply formula (\ref{e Om mth two}), which yields that the argument
 $\Arg \left( \Omega (\zeta_1,\zeta_2) - \om \right)$ of the QN eigenvalue's perturbation changes
 continuously for small
positive $\zeta_1$ and $\zeta_2$. This and (A2) imply that for small enough $\epsilon$ there
exist positive $\zeta_1$ and $\zeta_2$ such that $\zeta_1 + \zeta_2 = \epsilon$ and
 $\Arg \left( \Omega (\zeta_1,\zeta_2) - \om \right) $ equals $\theta$ (up to a multiple of $2 \pi$).
 Finally, we take such $\zeta_1$, $\zeta_2$ for a sequence $\epsilon_n$ going to $0$.
Note that the convexity of the families of admissible cavities and admissible directions is essentially
used here.

\section{Nonlinear eigenvalue problem for local extrema}
\label{s nonlin}

In this section we derive a nonlinear eigenvalue problem for QN modes corresponding
to local extremizers.

\vspace*{0.25cm}
\begin{lemma} \label{l xi}
Let $\ep$ be an admissible cavity and $\im z^2<0$. Then $\psi (s,z; \ep) \neq 0$ in the interval $0<s\leq l$,
 and the function $\xi (s) = \arg \psi^2 (s,z; \ep)$ is strictly
increasing (the values of $\arg$ are chosen such that $\xi$ is continuous for  $0<s\leq l$).
\end{lemma}

Let us show this. Assume first that $\psi (s_1) =0$, where $\psi (s) = \psi (s,z; \ep)$.
Then
\[
z^2 \int_0^{s_1} |\psi (s)|^2 \ep(s) \dd s =
- c^2 \int_0^{s_1} \psi'' (s) \overline{\psi (s)} \dd s
= - c^2  \int_0^{s_1} | \psi' (s) |^2 \dd s .
\]
The number in the right-hand side is real, but by our assumptions $z^2$ is not real.
Hence $s_1=0$ (note that $\psi' (0)=1$ and so $\psi$ is non-zero for small positive $s$).
Thus, $\psi (s) \neq 0$ for $s>0$. Note that
\[
\xi'(s) = 2 \im (\ln \psi (s))' = 2 \im \frac{\psi' (s) \overline{\psi (s)}}{|\psi(s)|^2} .
\]
This shows that, on the one hand, the function $|\psi|^2 \xi' $ tends to $0$ as $s$ goes to $0^+$,
and on the other hand,
\[
(|\psi (s) |^2 \xi' (s))' = 2 \im \psi'' (s) \overline{\psi (s)} =
- 2 c^{-2} \ep (s) |\psi (s)|^2 \im z^2 > 0.
\]
Hence $|\psi|^2 \xi'$ is positive for $0<s \leq l$, and so is $\xi'$.

Denote by
$\chi_{a,b}$ the function that equals $1$ in the interval $a <s< b $ and equals $0$ outside
this interval. For a finite number of complex numbers $z_1$, $z_2$, \dots, $z_n$, the convex hull of these
numbers, by definition, consists of numbers $c_1 z_1 + c_2 z_2 + \dots c_n z_n$, where
$c_1 + c_2 + \dots + c_n =1$ and all $c_j $ are nonnegative. Clearly, a set of complex numbers
is convex if and only if it contains a convex hull of every its finite subset.

Assume that $\ep$ is a piecewise continuous admissible cavity, and that
$\om$ with $\re \om >0$ is a QN eigenvalue of $\ep$.
Since changing of $\ep$ at a finite number of points does not influence QN eigenvalues,
this is not an essential restriction to assume from now on that $\ep$ takes
its limit from the left or its limit from the right at the points of discontinuity and the endpoints $s=0$ and $s=l$.

Further, suppose that $\ep_1<\ep (s)<\ep_2$
in a certain interval $s_1< s <s_2$. Then
\begin{equation} \label{e ep sep}
\ep_1+\delta_1<\ep (s)<\ep_2-\delta_1 \text{
in some narrower interval } s_3<s<s_4
\end{equation}
for small enough $\delta_1>0$.
Hence, the two directions
$\delta_1 \chi_{a,b} $ and $(-1) \delta_1 \chi_{a,b} $ are admissible for every $a$ and $b$ such that
$s_3 < a < b < s_4$.
It follows from Lemma \ref{l xi} that there exist $a_1$,$a_2$, $b_1$, $b_2$ with the properties:
\begin{itemize}
\item[(i)] $s_3 < a_1 < b_1 < a_2 < b_2 < s_4$,
\item[(ii)] the directional derivatives
$\pa_{\zeta} F (z; \ep + \zeta \chi_{a_1,b_1}) $, $ \pa_{\zeta} F (z; \ep + \zeta
\chi_{a_2,b_2} ) $ are non-zero and the difference of their arguments is not a multiple of $\pi$.
\end{itemize}
With such a choice of $a_{1,2}$ and $b_{1,2}$,
the convex hull of the four admissible directional derivatives
$\pa_{\zeta} F (z; \ep \pm \zeta \delta_1 \chi_{a_1,b_1}) $, $ \pa_{\zeta} F (z; \ep \pm \zeta
\delta_1 \chi_{a_2,b_2} ) $
contains a neighborhood of zero. Since the set of admissible directional derivatives is convex,
it also contains a neighborhood of zero. By Lemma \ref{l dir der} the pair $\{\om,\ep\}$ is not a
local extremizer.

Thus, we have shown that \emph{if an admissible pair $\{\om, \ep \}$
(with piece-wise continuous $\ep$) is a local extremizer, then
$\ep$ takes only the extreme possible values $\ep_1$ and $\ep_2$}.

Since the set of admissible directional derivatives at $\{\om,\ep\}$ is convex,
it either contains a neighborhood of zero
or, otherwise,
\begin{multline} \label{e in half-p}
\text{the set of nonzero admissible directional derivatives
is contained in } \\
\text{a certain complex half-plane } \ \theta_0 \le \arg z \le \theta_0+\pi.
\end{multline}

Assume again that $\{\om, \ep \}$ is a local extremizer and $\ep$ is piecewise continuous. Then Lemma \ref{l dir der} yields that the case (\ref{e in half-p}) takes place.
We have shown before that $\ep $ takes only values $\ep_1$ and $\ep_2$. Now, repeating the use of  Lemma \ref{l dir der},
we want to find restrictions on the points $s$ where
such an extreme $\ep$ switches between $\ep_1$ and $\ep_2$. To this purpose, consider directional derivatives produced
by directions $\pm \delta_1 \chi_{a,b}$ with intervals $a < s < b $ of small lengths.
 For any point $s_0$, there exists
an interval (of positive length) $s_1 \le s \le s_2$ containing $s_0$ and such that $\ep$ is constant on it.
To be specific, assume $\ep (s_0) = \ep_1$.
Then for any subinterval $a \le s \le b$ of the interval $s_1 \le s \le s_2$, the direction $\delta_1 \chi_{a,b}$ is admissible
whenever $0 \le \delta_1 \le \ep_2 - \ep_1$.
The corresponding admissible directional derivative
$\pa_\zeta F (\om, \ep + \zeta \delta_1 \chi_{a,b})$ equals
\[
\delta_1 \frac{\om^2}{c^2 \psi (l,\om;\ep)} \int_a^b  \psi^2 (s,\om;\ep) \dd s .
\]
Taking $a$ and $b$ close enough to $s_0$, it is possible to make the difference of arguments
\[
\arg \pa_\zeta F (\om, \ep + \zeta \delta_1 \chi_{a,b}) -
\arg \left( \frac{\om^2}{c^2 \psi (l,\om;\ep)} \psi^2 (s_0,\om;\ep)  \right)
\]
arbitrary small (after possible correction on a multiple of $2 \pi$).
Combining this with (\ref{e in half-p}), one can see that
\[
\theta_0 \le \arg \left( \frac{\om^2}{c^2 \psi (l,\om;\ep)} \psi^2 (s_0,\om;\ep)  \right) \le \theta_0+\pi .
\]
Similarly, if $\ep (s_0) = \ep_2$, then there exist $a$ and $b$ arbitrary close to $s_0$ such that
the direction $(-1) \delta_1 \chi_{a,b}$ is admissible whenever $0 \le \delta_1 \le \ep_2 - \ep_1$.
This yields
\[
\theta_0 - \pi \le \arg \left( \frac{\om^2}{c^2 \psi (l,\om;\ep)} \psi^2 (s_0,\om;\ep)  \right) \le \theta_0  .
\]

Summarizing, we see that, for a local extremizer $\{ \om, \ep \}$ with piece-wise continuous $\ep$, there exists
$\theta_1$ such that
\emph{$\ep$ takes the value $\ep_1$ in the intervals where $\theta_1 < \arg \psi^2 (s,\om;\ep) < \theta_1 + \pi$, and takes
the value $\ep_2$ in the intervals where $\theta_1 - \pi < \arg \psi^2 (s,\om;\ep) < \theta_1 $}
(here, the double inequalities with multi-valued $\arg$ function are assumed to be valid if they are valid for one of its values).
The values of $\ep$ at the finite set of points where $\psi^2 (s)$ cross the lines $\arg z = \theta_1 $ and $\arg z = \theta_1 + \pi$
are not important since these values do not influence QN eigenvalues of $\ep$.
The angle $\theta_1$ is related to $\theta_0$ of (\ref{e in half-p}) through
\begin{equation} \label{e th1}
\theta_1 = \theta_0 - \arg \frac{\om^2}{c^2 \psi (l,\om;\ep)}.
\end{equation}

Similar arguments also works for the case $\re \om < 0$ (for instance,
trough the use of the symmetry of QN eigenvalues w.r.t. $\ii \RR$) and
for the case $\re \om =0$
(for this case some modifications are needed, but considerations are simpler since the solution $\psi$
is real-valued).
Finally, note that the assumption of
piecewise continuity of $\ep$ can be dropped (that is, $\ep$ can be  assumed only to be an admissible $L^\infty$-function). This can be done if instead of small intervals $a < s < b$, one considers
sets of small measure and small diameter. Corresponding modifications require a number of purely technical details,
which a concerned reader can recover without difficulties using \cite{Ka12_pr}.
%(It occurs that this statement holds true also for general $L^\infty$-functions $\ep$
%from the family of admissible cavities.)

It is convenient to write the obtained result in terms of a nonlinear eigenvalue problem.
Denote by $\{ \im z > 0 \}$ the (open) upper complex half-plane, which, by definition, consists of complex numbers $z$ with positive imaginary part.
Similarly the half-plane $\{ \im z < 0 \}$ consists of numbers $z$ such that $\im z < 0$.
%and by $\overline{\{ \im z > 0 \}}$ the closed upper half-plane of $z$ such that $\im z \ge 0$.
Let $\ChiCpl  (\zeta) = 1$ when $\im \zeta > 0$, and $\ChiCpl  (\zeta) = 0$
when $\im \zeta \leq 0$.

\vspace*{0.25cm}
\begin{theorem} \label{t loc}
Assume that $\{ \om, \ep \}$ is a local extremizer (i.e., minimizer
or maximizer) of the decay rate for an admissible frequency
$\alpha= \re \om$. Then the boundary value problem (\ref{e bcE}) for the
nonlinear equation
\begin{equation} \label{e ep nl}
 E'' = - \frac{\om^2}{c^2} E \left[ \ep_1 + (\ep_2 - \ep_1) \ChiCpl
(E^2 ) \right]
\end{equation}
has a nontrivial solution $E(s)$ satisfying
\begin{equation} \label{e ep = chi}
\ep (s) = \ep_1 + (\ep_2 - \ep_1) \ChiCpl (E^2 (s) )  \ \ \ \text{ (almost everywhere) on } 0 < s < l.
\end{equation}
\end{theorem}

Let us show this for the case $\re \om >0$. Multiplying the solution $\psi$ to a constant $e^{\ii (\pi -\theta_1)/2}$ with
$\theta_1$ from (\ref{e th1}),
one gets another solution $E(s) = e^{\ii (\pi -\theta_1)/2} \psi (s,\om;\ep)$ to the linear equation
$ E'' = - \frac{\om^2}{c^2} \ep E $. This solution additionally satisfies (\ref{e bcE}) and  has the property that
$E^2$ takes values in $\{ \im z < 0 \}$ when $\ep (s)$ equals $\ep_1$, and takes values in $\{ \im z > 0 \}$
when $\ep (s) = \ep_2$
(after a possible correction of $\ep$ at a finite number of points).
So $\ep$ and $E$ are additionally connected by (\ref{e ep = chi}) and, therefore, $E$ is a solution to (\ref{e ep nl}).
Modifications of this proof for the cases $\re \om < 0$ and $\re \om =0$ are similar to that of \cite{Ka12_pr}.

Note that if $y(s)$ is a solution to the equation
\begin{equation} \label{e nonlin z}
 E'' = - \frac{z^2}{c^2} E \left[ \ep_1 + (\ep_2 - \ep_1) \ChiCpl
(E^2 ) \right],
\end{equation}
then $C y (s)$ is so for each \emph{positive} constant $C$.
This implies that any nontrivial solution $y$ to (\ref{e nonlin z}) satisfying $y(0)=0$ may be written in the form
$y(s) = C \Psi (s; z, \theta)$ with a positive constant $C$ and a function $\Psi (s; z, \theta)$ defined
by
\[
\pa_s^2 \Psi = - \frac{z^2}{c^2} \ \Psi \ \left[ \ep_1 + (\ep_2 - \ep_1) \ChiCpl
(\Psi^2 ) \right] , \ \
\Psi (0; z,\theta) = 0, \ \ \pa_s \Psi (0; z,\theta ) = e^{\ii \theta} ,
\]
where $z$ is a complex number  and $\theta$ is a number in the interval $-\pi <\theta \le \pi$.
It is easy to see that, for the initial data $E(0)$, $E'(0)$ such that $E(0)=0$ and $E'(0) \neq 0$, equation (\ref{e nonlin z}) has a unique solution, and so
$\Psi$ is well-defined.

For $-\pi <\theta \le \pi$ and a complex number $z$, define the function
\begin{equation*} \label{e W}
W (z,\theta ) = \ii z \frac{\sqrt{\ep_\infty}}{c} \Psi (l; z, \theta) -  \pa_s \Psi (l; z, \theta) .
\end{equation*}
One can see that the nonlinear eigenvalue  problem (\ref{e ep nl}), (\ref{e bcE}) has a nontrivial solution
exactly for those numbers $\om$ that satisfy
\begin{equation} \label{e W=0}
W  (\om, \theta ) = 0 \text{ for at least one } \theta .
\end{equation}

Let (\ref{e W=0}) be valid with a certain $\theta$. Then $y (s) = \Psi (s; \om, \theta)$ is a solution to
the original linear eigenvalue problem (\ref{e epE}), (\ref{e bcE})
with admissible $\ep (s) = \ep_1 + (\ep_2 - \ep_1) \ChiCpl
(y^2 ) $, and $\om$ is a QN eigenvalue of $\ep$. Thus, $\om$ is an admissible QN eigenvalue.
This, Theorem \ref{t loc}, and the definitions of optimal cavities and QN eigenvalues
(see Sect.\ref{ss def extr}) yield our main result.

\vspace*{0.25cm}
\begin{theorem} \label{t opt}
The minimal decay rate $\beta_{\min} (\alpha)$ for an admissible frequency $\alpha$
equals to the minimal real number $\beta$ having the property that
$W (\alpha - \ii \beta, \theta ) = 0$ for at least one $\theta$ in the range $-\pi <\theta \le \pi$.

Moreover, let $\om = \alpha - \ii \beta_{\min} (\alpha)$ be the corresponding optimal QN eigenvalue.
Then for any $\theta$ such that $W (\om , \theta ) = 0$,
the cavity defined by
\begin{equation} \label{e ep opt =}
\ep (s) = \ep_1 + (\ep_2 - \ep_1) \ChiCpl \left( \Psi^2 (s; \om, \theta) \right)
\end{equation}
is optimal for the frequency $\alpha$. And vice versa, for each cavity $\ep$ optimal for $\alpha$ ,
there exists $\theta$ such that (\ref{e ep opt =}) holds.
\end{theorem}

\section{Conclusions}
The paper is concerned with QN eigenvalues of 1-D leaking optical
cavities. The cavities are described by dielectric permittivity
function $\ep$ which depends on one variable $x_3$ and is assumed to
take values in a fixed range $\ep_1 \le \ep (x_3) \le \ep_2$. Cavities
satisfying this restriction are called admissible.

We study analytically cavities that produce QN eigenvalues $\om$
with a given frequency $\re \om = \alpha$ and locally extremal (locally maximal or minimal) for this frequency
decay rate $\beta = - \im \om$. We show that cavities with such locally extremal properties are 1-D photonic crystals
consisting of alternating layers of two materials with two extreme
allowed dielectric permittivities $\ep_1$ and $\ep_2$.
This explains effects observed in numerical experiments for very kindred (but slightly different)
optimization problems \cite{KS08,HBKW08}.
To find thicknesses of the layers in extremal cavities,
we derive a nonlinear eigenvalue problem for their QN modes.
It occurs that $x_3$ coordinates of interface planes between the layers are tied
to rotation in the complex plane (i.e., to the $\arg$-function) of the corresponding extremal QN mode.

For each admissible frequency $\alpha$ there exist a cavity that creates a QN
eigenvalue with minimal possible decay rate.
In the paper, this QN eigenvalue is called optimal for the frequency $\alpha$.
We show that such optimal QN eigenvalues can be easily found via zeroes of a function $W$,
which is constructed by solutions of the above mentioned nonlinear equation
and depends only on two variables. This effectively excludes the unknown optimal dielectric permittivity $\ep$
from the process of calculation of an optimal QN eigenvalue. After calculation of an optimal QN eigenvalue and
corresponding QN modes,
optimal dielectric permittivity functions can be easily recovered from their connection with optimal QN modes.

\vspace*{0.5cm}
\section*{Acknowledgements}
The author is grateful to Igor Chueshov, Robert Kohn, and Roald Trigub
for interest to this research and stimulating discussions, and to Andrey Shishkov for constant support.

%*******************************************************************
%BIBLIOGRAPHY
%*******************************************************************

\end{document}